\documentclass[12pt]{article}
\usepackage{amsmath}
\usepackage{amssymb}
\usepackage{amsthm}
\usepackage{amsfonts}
\usepackage{graphicx}

\usepackage{bbm,epsfig,graphics,epic,color,rotating,color}
\numberwithin{equation}{section}

\newcommand{\R}{{\mathbb R}}

\newcommand{\x}{{\mathbf x}}

\newcommand{\wt}[1]{\widetilde{#1}}
\newcommand{\wh}[1]{\widehat{#1}}

\newcommand{\ed}{\mathrm{d}}

\newcommand{\reff}[1]{(\ref{#1})}

\newtheorem{theorem}{Theorem}[section]

\theoremstyle{definition}

\title{Dynamics of Tectonic Plates}

\author{E. Pechersky$^{1}$ S. Pirogov$^{1}$,  G. Sadowski$^2$ and A. Yambartsev~$^{3}$    }

\begin{document}

\maketitle{\footnotesize
\noindent $^{1}$ Dobrushin laboratory of Institute for Information
Transmission Problems of Russian Academy of Sciences,\\
19, Bolshoj Karetny, Moscow, Russia.\\
E-mail: pech@iitp.ru, pirogov@iitp.ru

\noindent $^2$ Department of Mineralogy and Geotectonics, Institute of Geosciences
University of S\~ao Paulo\\
Rua Prof. Guilherme Milward 246
S\~ao Paulo, SP, Brazil,  05506-000\\
E-mail: sadowski@usp.br

\noindent $^3$ Department of Statistics, Institute of Mathematics
and Statistics, University of S\~ao Paulo\\ Rua do Mat\~ao 1010,
CEP 05508--090, S\~ao Paulo SP, Brazil.\\
E-mail: yambar@ime.usp.br }

\begin{abstract}
We suggest a model that describes a mutual dynamics of tectonic plates. It is a sort of stick-slip dynamics modeled by a Markov random process. The process defines the dynamics on a micro-level. A macro-level is obtained by a scaling limit which leads to a  system of integro-differential equations which determines a kind of mean field systems. The conditions when the Gutenberg-Richter empirical law holds are presented at the mean field level. Those conditions are rather universal and independent of the features of the resistant forces.
\end{abstract}

\section{Introduction.}
Tectonic plate construction of the earth lithosphere is generally accepted  and recently well described (for example, see \cite{KC,Pr}). Moreover, the tectonic plates of another planets of the solar system, such as Mars and Venus, are a subject of investigations as well \cite{Sup}.  The tectonic plate motion as a reason of earthquakes is also a widely held view.  However, the mechanism of the earthquakes emergence in a course of the plate motions is subject of intensive investigations in Geosciences and Geophysics as well as in Mathematics and Physics (see \cite{BK, CL, SBB, PP,  SS}).

The most popular models represent a \textit{stick-slip} motion.
Brace and Byerlee \cite{BB} suggested that earthquakes must be the result of a stick-slip frictional instability rather then caused by fracture appearance and propagation. The  earthquake being the result of a sudden slip along a pre-existing fault or plate surface and the ``stick" being  the interseismic period of  strain accumulation.  For example San Andreas Fault or Benioff-Wadati subduction zones, seem to exhibit sudden slip  events followed by ``silent" slip and renewal of the contact points creating a stick behavior. 

This idea has been applied in several laboratory experiments which lead to models consistent with  Ruina-Dietrich ``rate and state-variable friction law" \cite{Sh}, \cite{Ruina}. Another theoretical model following this idea was proposed in \cite{BK} (Burridge-Knopoff model), where experimental and theoretical results are discussed. The model in \cite{BK} is a one-dimension  chain of massive blocks tied by springs and situated on a rough unmoving surface.  The heading block is pulled with a constant speed. Then the train of rest blocks shows the stick-slip behavior under some conditions related to resistance forces between blocks and the rough surface.

In this paper we suggest a model for the stick-slip dynamics. The model describes a stochastic dynamics on a micro-level. It is the stochastic dynamics of a point field which is a set of contact points on surfaces.  A combination of three types of changes forms the dynamics of the plate: 1) a deterministic motion causing a deformation of the contact points and the increase of strain; 2) an appearance of new contacts that changes the dynamics; 3) break of the existent contacts which causes a jump-wise change of the dynamics velocity. Some of physical phenomena are described in the model by  a number of  parameters. One of the parameters is a rate of contact destructions. We show that these parameters are  related to the Gutenberg-Richter law (see \cite{Ut}). For the elastic resistance we assume a linear dependence of the force on a value of the contact deformation.  The model defined and numerically analyzed in \cite{Braun} is the most similar model we found in the literature.

The micro-level model defined by a Markov piece-wise deterministic random process is a base for a macro-level model represented by a system of integro-differential equations. This model is rather universal to describe tectonic plate motions as well as the friction between different matter plates.  The  mean field method used in this paper was also applied in \cite{KPP}.

\vspace{1cm}

Sections \ref{s2} and \ref{s3} contain the descriptions and the basic results including the condition when the Gutenberg-Richter law is satisfied. This condition requires a special asymptotic law for the breaks of the contact points which we call \textit{inverse deformation law}. Section \ref{s4} contains some rigorous descriptions.

\section{The model description. Finite volume.}\label{s2}

We consider a plate moving on a solid unmoving substrate with contact area $\Lambda$. The plate is subjected to the action of a constant force $F$ in some fixed direction  (a \textit{moving force}). The moving force causes a motion of the plate in the direction of $F$. Asperities of both plate and solid substrate can be contacted, what creates a resistance force due to deformation of the contacting asperities. The set of {\it contact points} we denote by $\Omega$ and its elements by $\omega$. The set $\Omega$ is a finite subset of $\Lambda$.

\subsubsection*{Deterministic part} In the course of the  plate motion, any contacting asperity is being deformed until its failure. 

Any contact point $\omega \in \Omega$ creates a resistant force  depending on magnitude of the asperity deformation. In a general setting the resistant force is an expression $\varkappa \min\{x_\omega,x_\omega^{\alpha}\}$, where is a constant from $[0,1]$.   When $\alpha=1$ then the resistance is purely elastic, and $\varkappa$ is Hooke's constant. The sum of all  resistance forces over all contact points gives  the total resistance force. Therefore the resultant force acting on the plate is 
\begin{equation}\label{force}
G=\Bigl[F-\varkappa\sum_{\omega \in \Omega} \min\{x_\omega,x_\omega^{\alpha}\} \Bigr]_+,
\end{equation}
where $[A]_+=\max\{A,0\}$. This means that the resistance force cannot be greater than the moving force.

The force $G$ causes a deterministic motion of the plate. We assume that the plate is a monolithic hard object any point of which is moving with the same velocity $v$. The next our assumption implies that  the dynamics of the plate under the force action follows  the so called Aristotle mechanics. It means that the velocity $v$ of the motion is proportional to the acting force,
\begin{equation}\label{Arist}
v=\gamma G,
\end{equation}
in contrast to Newton mechanics where the acceleration is proportional to the acting force. Here $\gamma$ is constant. This assumption is due to the fact that the moving plate is immersed in a viscous medium and therefore inertia, despite the large mass of the plate, has little influence on the nature of the motion. 

If  during a time interval $[t_1,t_2]$ the set $\Omega$ is not changed (a new contact does not appear and any existing contact does not disappear) then we have a deterministic dynamics on $[t_1,t_2]$ evolving according to the equation
\begin{equation}\label{dyn2}
\frac{\ed x_\omega(t)}{\ed t}=v(t)=\gamma\Bigl[F-\varkappa\sum_{\omega' \in\Omega} z_{\omega'}^{(\alpha)}(t)\Bigr]_+
\end{equation}
with an initial velocity value $v(t_1)$, (see \reff{force} and \reff{Arist}). Here we denote $z_\omega^{(\alpha)}=\min\{x_\omega,x_\omega^{\alpha}\}$.

In the case of the elastic resistance forces ($\alpha=1$) , if 
\begin{equation}\label{ela}
v(t_1)=\gamma\Bigl[F-\varkappa\sum_{\omega' \in\Omega} x_{\omega'}(t_1)\Bigr]_+=\gamma\Bigl(F-\varkappa\sum_{\omega' \in\Omega} x_{\omega'}(t_1)\Bigr)>0
\end{equation}
then a solution of \reff{dyn2}  on $[t_1,t_2]$ is
\begin{equation}\label{vel}
v(t)=v(t_1)e^{-\varkappa \gamma |\Omega|(t-t_1)}
\end{equation}
for $t\in[t_1,t_2]$, where $|\Omega|$ is a number of the points (the contacts) in the set $\Omega$ (see rigorous evaluations in the section \ref{s41}).

This dynamics determines the motion of the plate when the set of contact points is fixed. 

\subsubsection*{Stochastic part}The deterministic dynamics is interrupted by random events of two kinds: either some contact from $\Omega$ disappears or a new contact appears. In both cases the number of the contacts is changed from $n=|\Omega|$ to either $n-1$ or $n+1$.
Next we describe the dynamics of  the appearance and disappearance (birth and death) of the contact points. These dynamics are random and of Markov type.

Assume that there were no any random events on the interval $[t_1,t_2]$, and at $t_2$ a new contact appears.  Its deformation is equal to 0 at the moment $t_2$. Let $\Delta_\Omega=\{x_\omega\}$ be a set of all contact point deformations. 
The contact set $\Omega$ is changed to a set $\Omega'$ at the moment $t_2$, and the new deformation set is $\Delta_{\Omega'}=\Delta_\Omega\cup\{0\}=\{0,x_\omega:\: \omega\in\Omega\}$.

It is assumed that the birth of new contact points at $t_2$ depends on the velocity $v(t_2)=\gamma[F-\varkappa\sum_{\omega' \in\Omega} z_{\omega'}^{(\alpha)}(t_2)]_+$ of the plate. This dependence is determined  by a birth 
rate  ${c}_b(v(t))>0$.
The function ${c}_b$ should reflect  physical properties of the plate such as a fractal dimension of the asperities and many others.  It is clear that new contacts do not appear if the velocity $v$ is zero, hence ${c}_b(0)=0$. A natural choice of ${c}_b$ is ${c}_b(v)=\bar c_b v$, linear dependence on the velocity, where $\bar c_b>0$. 

Because the new contact deformation is 0 then 
$$\sum_{\omega' \in\Omega} z_{\omega'}^{(\alpha)}(t_2)=\sum_{\omega' \in\Omega'} z_{\omega'}^{(\alpha)}(t_2).$$ 
Thus the velocity $v(t)$ is continuous at $t=t_2$. However the velocity derivative $\frac{\ed v(t)}{\ed t}\big|_{t=t_2}$ at $t_2$ is discontinuous; see Figure~\ref{F1} of a typical path of the velocity. The left derivative at $t_2$ is greater than the right one at the same time $t_2$. In the case of elastic resistance forces:
\begin{eqnarray*}
\lim_{t\uparrow t_2}\frac{\ed v(t)}{\ed t} &=& -\gamma\varkappa nv(t_1)e^{-\gamma\varkappa n(t_2-t_1)} \\
&>& \lim_{t\downarrow t_2}\frac{\ed v(t)}{\ed t}=-\gamma\varkappa (n+1)v(t_1)e^{-\gamma\varkappa n(t_2-t_1)},
\end{eqnarray*}
where $n=|\Omega|$.

The disappearance (death) of a contact $\omega$ from $\Omega$ is determined by the rate $c_u \equiv c_u( x_\omega)$, it may depend on the deformation $x_\omega$. If at the moment $t_2$ the contact $\omega$ disappears then the velocity has discontinuity at $t_2$. The velocity increases abruptly by the value $\gamma\varkappa z_\omega^{(\alpha)}$. We shall assume that the velocity is continuous from the right at $t_2$, that is  
$$
\lim_{t\downarrow t_2} v(t)=v(t_2).
$$

\subsubsection*{Complete view}We describe now a complete evolution of the plate velocity by the dynamic described above. The time is split into the intervals $[0,\infty]=\bigcup[t_i,t_{i+1}]$ such that during every interval $[t_i,t_{i+1}]$ the plate is moving deterministically according to \reff{dyn2} and \reff{vel} substituting $t_1$ by $t_i$. The set of the points $R=\{t_i\}$ is the moments of random events: either a new contact appears or one of the existing contacts disappears. The set of the random moments $R$ is splitted by the set $R_b$ of the appearing contacts and the set $R_u$ of the disappearing contacts, $R=R_b\cup R_u$. 

Any disappearance of the contact releases an energy which switches to an energy of seismic waves in the plates. The amount of released energy depends on the deformation value $x_\omega$ of the disappeared contact. We assume that the seismic wave  energy is proportional to the primitive of the resistance force (the potential). For large deformation $x$ it is proportional to $x^{1+\alpha}$ . If $x_\omega$ is large then the oscillation amplitude may be large which can be observed as an earthquake.

\begin{figure}%\label{F1}
\begin{center}
\includegraphics[scale=0.5]{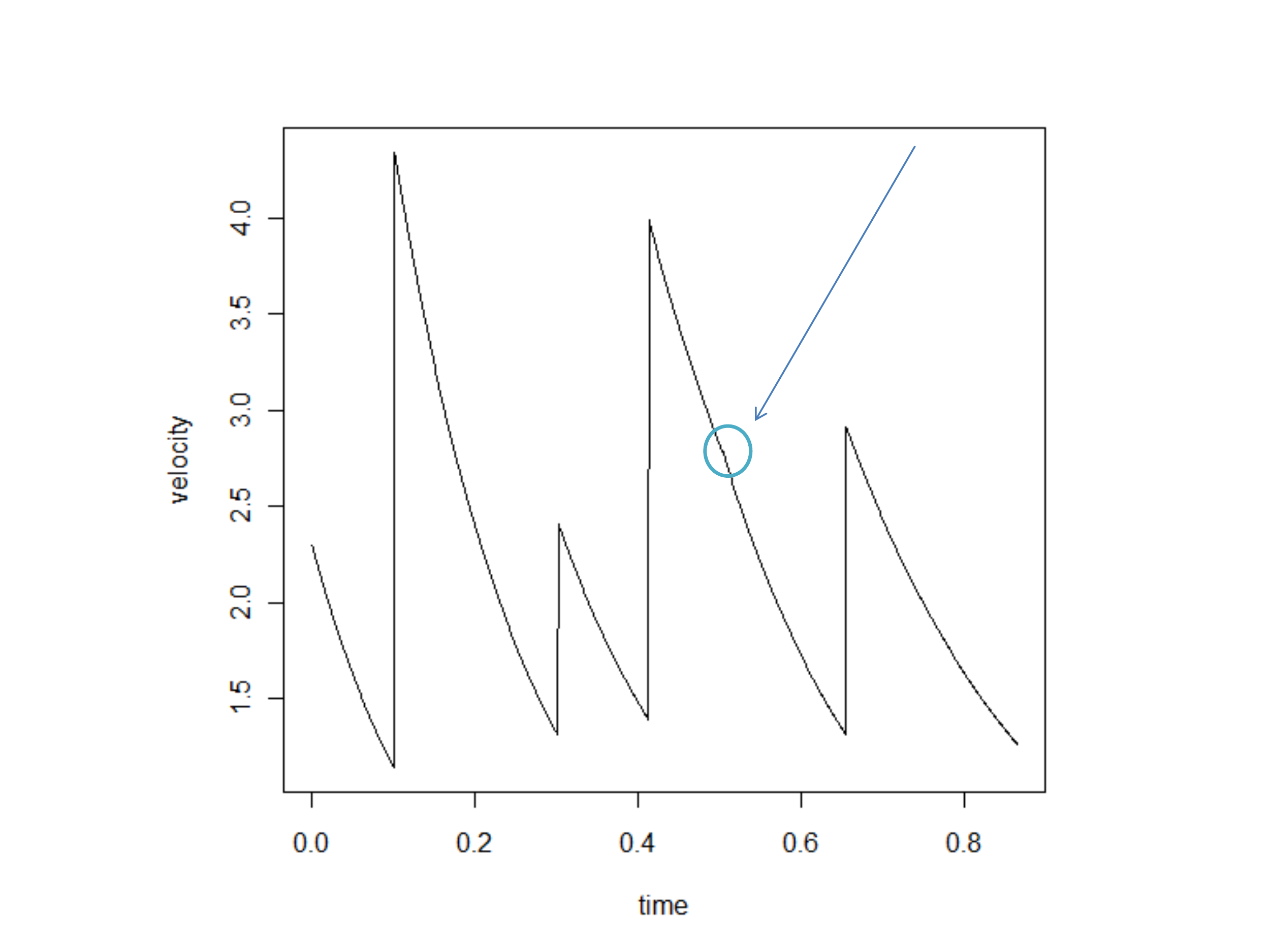}
\end{center}
\caption{A typical path of the velocity $v(t)$ is a piece-wise continuous function between two nearest
deaths of the contacts. The velocity has  jumps at  moments of the death of particle. The derivative of the 
velocity increases abruptly at the moment of the birth of a particle (blue circle). }%
\label{F1}
\end{figure}

A typical random path of the velocity is presented on Figure \ref{F1}.

%\newpage

\section{Scaling limit. Infinite volume}\label{s3} In this section we propose an analytical approach which allows to study some properties of the defined model. The idea is to consider a very large number of the asperities. It means that we  consider a limit  of the size of the contact area $\Lambda$ going to infinity.  In the limit we obtain infinite number of asperities and then a distribution of the asperity deformations  is described by a density function $\rho(x,t)$, $\rho(x,t)\ed x$ has a  meaning of a  number  of the asperities having their deformation in the interval $\ed x$  at the moment $t$. To obtain a  reasonable model in this limit we have to change the values of the parameters determining the plate motion model. Some of the parameters must depend on the size of the contact surface $\Lambda$. Namely, the birth intensity $c'_b= \bar{c}_b\Lambda$ and the acting force $F'=F\Lambda$, are proportional to the size $\Lambda$. The death parameter we take without changes $c'_u=c_u$. Further we omit the sign ${}'$.

A system of  equations describing a behavior of the density $\rho(x,t)$ and the velocity $v(t)$ derived from balance considerations is
\begin{equation}\label{y-3}
    \frac{\partial\rho(x,t)}{\partial t} + {v(t)} \frac{\partial\rho(x,t)}{\partial x} = {\bar c}_b {v(t)} \delta_0(x) - {c}_u(x)\rho(x,t),
\end{equation}
and
\begin{equation}\label{y-4}
v(t)= \gamma \left( F - \varkappa \int_{-\infty}^{+\infty} z^{(\alpha)}(x)\rho(x,t)dx \right)_{+},
\end{equation}
where $z^{(\alpha)}(x)=\min\{x,x^{\alpha}\}$.
Here we define the function $\rho(x,t)$ on all line $x\in\R$ such that it is equal to 0 when $x<0$.  

In the two next subsections we study some properties of the solutions of \reff{y-3}, \reff{y-4}.

\subsection{Stationary state and the Gutenberg-Richter law} In this subsection we study some properties of \reff{y-3} in a stationary regime, when $\rho(x,t)$ does not depend on time $t$. The equation \reff{y-3} is then
\begin{equation}\label{y-3.1}
    v \frac{\ed\rho}{\ed x}(x) = \bar c_b v \delta(x) -  c_u(x) \rho(x),
\end{equation}
(cf \reff{eqsta}).
Because of the stationarity the velocity $v$ of the plate does not depend on time $t$. It means that $v$ is a constant in \reff{y-3.1}.  
The solution of \reff{y-3.1} connects $\rho(x)$ and $v$ as the following
\begin{equation}\label{y-12}
    \rho(x) = \left\{ \begin{array}{ll} {\bar c_b} \exp \bigl\{ -\frac{1}{v} \int_0^x c_u(y) dy \bigr\}, & \mbox{ when }x>0, \\
    0, & \mbox{ when }x<0. \end{array} \right.
\end{equation}
The density at $x=0$ is $\rho(0)=\bar c_b$. For $v=0$ the stationary regime is trivial, $c_u(x)\rho(x)=0$.

When an asperity subjected to a  deformation size $x$ disrupts then it releases the energy
\begin{equation}\label{3.5}
e(x)=\mu x^{1+\alpha}x
\end{equation}
for large $x$, where  $\mu$ a constant of proportionality. 

\vspace{.5cm}

\textit{The Gutenberg-Richter law}. The Gutenber-Richter law is observed, in practice, on a restricted area of energy values $[e_1,e_2]$. Therefore we study the intensity distortion  function $c_u(x)$  on  corresponding the deformation interval $[x_1,x_2]$, where $x_i=e_i^{\frac1{1+\alpha}}$. It is assumed that $x_1>1$. Let $n(e)\ed e$ be a distribution of a number of the distortions along the  energy $e$ axis, and let $m(x)\ed x$ be a distribution of number of the distortions arising from the asperities destroyed at a size $x$ of their deformations. The Gutenberg-Richter law claims that 
\begin{equation}
n(e)\ed e\propto \frac{1}{e^w}\ed e,
\end{equation}
where $w>1$. Observations show that $w$ lies in the range $1.7-2.1$ (see \cite{Ut}). Changing variables defined by \reff{3.5} we obtain 
\begin{equation}\label{3.7}
n(e(x))\ed e(x)\propto n(e(x))x^{\alpha}\ed x\propto  \frac{1}{e^w(x)}x^{\alpha}\ed x
\end{equation}
The distribution $m(x)\ed x$ can be expressed by the density $\rho(x)$ of the deformations and the intensity $c_u(x)$ of the distortions  as
\begin{equation}
m(x)\ed x\propto c_u(x)\rho(x)\ed x.
\end{equation}
Remark next that 
\begin{equation}
m(x)\ed x=n(e(x))x^{\alpha}\ed x.
\end{equation}
It follows now from \reff{3.7} that 
\begin{equation}\label{3.10}
c_u(x)\rho(x)\ed x\propto \frac{1}{e^w(x)}x^{\alpha}\ed x=\frac{x^{\alpha}}{x^{w(1+\alpha)}}\ed x.
\end{equation}

The left hand side of this relation is the frequency (density) of the demolished contacts achieved with deformation values equal at $\ed x$. 

The equation   \reff{3.10} shows a mutual behavior dependence of  $c_u(x)$ and $\rho(x)$ for large $x$ (at least, $x>1$):
\begin{equation}\label{ddens}
\rho(x)\propto\frac{1}{x^pc_u(x)},
\end{equation}
where  $p=w(1+\alpha)-\alpha$. The next analysis uses the function $c_u(x)$ on the interval $[0,x_2]$. The function $c_u(x)$ on $x\in[0,x_1]$ is assumed arbitrary continuous one. The relation 
\begin{equation}\label{y-13}
\bar c_b\exp \left\{ -\frac{1}{v} \int_0^x c_u(y) \ed y \right\}\propto\frac{1}{x^pc_u(x)}
\end{equation}
for $x\in[x_1,x_2]$ follows from \reff{y-12}. That is
\begin{equation}\label{y-13.1}
\bar c_b\exp \left\{ -\frac{1}{v} \int_0^x c_u(y) \ed y \right\}=\frac{A}{x^pc_u(x)}
\end{equation}
for some $A>0$.
Taking logarithm and derivative of \reff{y-13.1} we obtain the differential equation
\begin{equation}\label{y-14}
\frac1{v}c_u(x)=\frac{p}{x}+\frac{c'_u(x)}{c_u(x)}
\end{equation}
on the interval $[x_1,x_2]$, where $c_u'(x)$ means the derivative of $c_u(x)$. 
We assume an initial value $c_u(x_1)>0$ is given. 

\subsubsection*{The rate of the disruptions corresponding to Ritenberg-Richter law}
\textit{Assume that the Gutenberg-Richter law is satisfied on $[x_1,x_2]$. Under this assumption  the general solution of \reff{y-14} is
\begin{equation}\label{idl}
c_u(x)=\frac{(p-1)v}{x+Bx^p},
\end{equation}
where $B\geq 0$.}

The fact that \reff{idl} is the solution of \reff{y-14} can be verified directly. 
%\subsubsection*{Inverse deformation law} 

The relation \reff{idl} shows that the disrupting asperity rate on $[x_1,x_2]$ is decreasing: \textit{greater the asperity deformation $x$ less a probability to be disrupted}. On other hand this correspondence  means: greater the asperity deformation slower the density $\rho(x)$ decreasing. In the case $B>0$ the density $\rho(x)$ can not be done as small as possible since $\lim_{x\to\infty}\rho(x)>0$ (see \reff{ddens} and \reff{idl}. Such $\rho$ has no a physical sense, thus we consider the case $B=0$. 

\subsubsection*{The inverse deformation law.} 

\textit{If $B=0$ then the disruption rate $c_u$ is decreasing hyperbolically, }
\begin{equation}\label{16}
c_u(x)=\frac{(p-1)v}{x}
\end{equation}

Assume that $c_u(x)$ is defined on $[0,\infty)$ such that \reff{16} holds on $[x_1,x_2]$ the following equations define $v$ and $\rho(x)$:
\begin{eqnarray}\label{3.16}
v&=& \gamma \left( F - \varkappa \int_{-\infty}^{+\infty} z^{(\alpha)}(x)\rho(x)dx \right)_{+},\\
\vphantom{\int}\rho(x) &=& {\bar c_b} \exp \bigl\{ -\frac{1}{v} \int_0^x c_u(y) dy \bigr\},  \mbox{ for }x>0. \label{3.18}
\end{eqnarray}
If solutions of \reff{3.16} and \reff{3.18} are such that $v>0$ then
\begin{equation}
\rho(x)=
\begin{cases}
{\bar c_b}\exp\left\{-\frac1{v}\int_0^xc_u(y)\ed y\right\}&\mbox{ if }x<x_1,\\
{\bar c_b}\exp\left\{-\frac1{v}\int_0^{x_1}c_u(y)\ed y\right\}\left(\frac{x_1}{x}\right)^{p-1}&\mbox{ if }x_1\leq x<x_2,\\
{\bar c_b}\exp\left\{-\frac1{v}\int_0^{x_1}c_u(y)\ed y\right\}\left(\frac{x_1}{x_2}\right)^{p-1}\exp\left\{-\frac1{v}\int_{x_2}^{x}c_u(y)\ed y\right\}&\mbox{ if }x_2\leq x.
\end{cases}
\end{equation}
and
\begin{eqnarray*}
&v=&\gamma \left( F - \varkappa{\bar c_b}\left[\int_0^{x_1}z^{(\alpha)}\exp\left\{-\frac1v \int_0^xc_u(y)\ed y\right\}\ed x\right.\right.+\\
&&\frac{x_1^{p-1}}{p-\alpha}\exp\left\{-\frac1v \int_0^{x_1}c_u(y)\ed y\right\} \left[x^{\alpha-p}_1- x^{\alpha-p}_2\right] +\\
&&\left.\left.\left(\frac{x_1}{x_2}\right)^p\int_{x_2}^\infty x^{\alpha}\exp\left\{-\frac1v \int_{x_2}^{x}c_u(y)\ed y\right\} \ed x \right]\right).
\end{eqnarray*}
Recall that $p=w(1+\alpha)-\alpha$

\vspace{1cm}
 
\subsubsection{Examples} 
\textit{Example 1.} Assume that $x_1=1$ and the constants $F,\ \gamma,\ \varkappa,\ {\bar c_b},\  f(\alpha)$  and $x_2$ are such that 
\begin{equation}\label{3.21}
\wt v=\gamma \left( F - \varkappa{\bar c_b}\left[\frac{1}{1+\alpha}+\frac1{p-\alpha} \left[1-\frac 1{x^{\alpha-p}_2}\right]  \right]\right)>0.
\end{equation}
We consider the following rate of the contact destruction
\begin{equation}\label{ex1.1}
c_u(x)=\begin{cases}0,&\mbox{ if }x\leq 1,\\
\frac{\wt v(p-1)}{x},&\mbox{ if }1<x\leq x_2\\
\infty,&\mbox{ if }x>x_2.
\end{cases}
\end{equation}
Then the solutions of  \reff{3.16} and \reff{3.18} are $v=\wt v$ and
\begin{equation}
\rho(x)=
\begin{cases}
{\bar c_b}&\mbox{ if }x<1,\\
{\bar c_b}\left(\frac{1}x\right)^{p-1}&\mbox{ if }1\leq x<x_2,\\
0&\mbox{ if }x_2\leq x. 
\end{cases}
\end{equation}

\textit{Example 2.} Assume that $x_1=1$ and $\alpha=1$. Let the constants $F,\ \gamma,\ \varkappa,\ {\bar c_b}$  and $x_2$ be such that there exists a constant $\wt a>0$ such that 
\begin{equation}
F>\wt r=\bar c_b\varkappa\left[\frac1{\wt a^2}\big(1-(1+\wt a)e^{-\wt a}\big)+\frac{e^{-\wt a}}{p-1}\big(1-x_2^{1-p}\big)\right].
\end{equation}
Let
\begin{equation}
c_u(x)=\begin{cases}
\wt a\gamma(F-\wt r),&\mbox{ if }x<1\\
\frac{(p-1)\gamma(F-\wt r)}{x},&\mbox{ if }1\leq x<x_2,\\
\infty,&\mbox{ if }x_2\leq x
\end{cases}
\end{equation}
then the solutions of \reff{3.16} and \reff{3.18} are
\begin{eqnarray*}
&&v=\gamma(F-\wt r)\\
&&\rho(x)=\begin{cases}
e^{-\wt ax},&\mbox{ if }x<1,\\
e^{-\wt a}\left(\frac{1}{x}\right)^{p-1},&\mbox{ if }1\leq x<x_2,\\
0,&\mbox{ if }x_2\leq x.
\end{cases}
\end{eqnarray*}
For the considered case $p=2w-1$.

There existence of the constant $\wt a$ follows from a limit
\begin{equation}
\lim_{a\to\infty}\frac1{ a^2}\big(1-(1+a)e^{- a}\big)+\frac{e^{-a}}{p-1}\big(1-x_2^{1-p}\big)=\infty
\end{equation}

%%%%%%%%%

\medskip

\subsection{Non-stationarity}Finding a solution of \reff{y-3} and \reff{y-4} in general case is rather a difficult problem. We describe a non-stationary behavior of the system in a particular case when the disruption intensity, $c_u(x)\equiv \bar c_u$ does not depend on the deformation value $x$. This model was studied in \cite{PSSY} on the micro-level. Here we obtain the same results on the macro-level.

We reduce the study of \reff{y-3} and \reff{y-4} to a dynamical system which can be completely investigated.
For formal considerations  see in Section \ref{s4}

\subsubsection{A dynamical system}

Define $$
N(t)=\int \rho(x,t)dx \  \  \mbox{ and } \  \  M(t) = \int x \rho(x,t) dx.
$$
Then $v(t)=\gamma \bigl( \bar F - \kappa M(t) \bigr)_+$ (see \reff{y-4}). 

We assume and check later that $\rho(x)\to 0$ and $x\rho(x)\to 0$, when $x\to \infty$.

Integrating (\ref{y-3}) over $x$   we obtain
\begin{equation}\label{y-5}
    \frac{dN(t)}{dt} = \bar c_b v(t) - \bar c_u N(t)
\end{equation}
Multiplying (\ref{y-3}) by $x$ and integrating, we obtain
\begin{equation}\label{y-6}
    \frac{dM(t)}{dt} +  v(t) \int x \frac{\partial\rho(x,t)}{\partial x} dx = - \bar c_u M(t),
\end{equation}
because of $\int x\delta(x)\ed x =0$. Evaluating the integral in the left-hand side of (\ref{y-6}) by parts we obtain the following system of equations on the plane $(N,M)$:
\begin{equation}\label{y-7}
    \left\{
    \begin{array}{rcl}
    dN/dt &=& \bar c_b v - \bar c_u N \\
    dM/dt &=&   v N - \bar c_u M
    \end{array} \right.
\end{equation}
describing a dynamical system in the plane. 

In the quarter-plane $(M,N),~M>0,~N>0$, there exists a unique point $(M_0,N_0)$ which is a fixed point of \reff{y-7}. It means that $(M_0,N_0)$ is a solution of  \reff{y-7} at the assumption that $\frac{\ed M}{\ed t}=\frac{\ed N}{\ed t}=0$. This solution is 
\begin{eqnarray}\label{stat}
% \nonumber to remove numbering (before each equation)
  M_0 &=& \frac{F}{\kappa}- \frac{F}{2\gamma\kappa a}\left[\sqrt{1+4\gamma a}-1\right]  \\
  N_0 &=& \frac{\bar c_b}{\bar c_u}\frac{F}{2a} \left[\sqrt{1+4\gamma a}-1\right],
 \nonumber
\end{eqnarray}
where 
\begin{equation}\label{orde}
a=\gamma\kappa F\frac{\bar c_b}{\bar c_u^2}.
\end{equation}
The fixed point $(M_0,N_0)$ is stable, that is the dynamic \reff{y-7} is such that a point $(M(t),N(t))$ is attracted to $(M_0,N_0)$ if $(M(t),N(t))$ is in a neighborhood of $(M_0,N_0)$. There are two different ways how a path $(M(t),N(t))$ is moving to $(M_0,N_0)$. The type of the ways depends on a value of $a$. This behavior is investigated by a linearization of the non-linear equations \reff{y-7} (see Section \ref{4.3}). At the fixed point, the density $\rho(x,t)$ does not depend on time and is equal to
\begin{equation}\label{stat2}
\rho(x)=\bar c_b\exp\left\{-\frac{\bar c_u}{v_0}x\right\},
\end{equation}
being a solution of the stationary version 
\begin{equation}\label{eqsta}
 \frac{\partial\rho(x)}{\partial x} = {\bar c}_b {v_0} \delta(x) - {\bar c}_u\rho(x)
\end{equation}
of \reff{y-3}. The velocity $v_0$ at the fixed point is a solution of the cubic equation
$$
v=\gamma F- \gamma \kappa\frac{\bar c_b}{\bar c_b^2}v^3,
$$
which has an unique positive solution between $v=0$ and $v=\sqrt[3]{\frac {F\bar c_u^2}{\varkappa \bar c_b}}$ (see \reff{y-4}).

\vspace{.3cm}

Remark that in the stationary regime  $\rho$ is exponentially decreasing (see \reff{stat2}).  Therefore Gutenberg-Richter law is not satisfied for a constant intensity $\bar c_u$ of the disruptions.

\newpage

\section{Mathematical tools}\label{s4}In this section we present  some mathematical justifications of the facts described in the previous sections.

\subsection{Markov process}

The model description, see Section \ref{s2}, shows that the stochastic dynamics of the plate is a Markov process which we define hereby.

Remark that the model does not care about the positions of the points $\omega$ in the area $\Lambda$, but it essentially depends on the displacements $x_\omega$ of any contact point $\omega$ (see \reff{force} and \reff{dyn2}). Thus a state of the moving plate is described by a set
\begin{equation}\label{conf}
\mathbf{x} := \{ x_\omega\}_{\omega \in \Omega} \subset \R_+=[0,\infty).
\end{equation}
of all displacements of the contacts $\Omega$.

It is clear then that we consider one-dimensional model. 

Further we will omit the index $\omega$. Let $X$ be a set of all finite configurations $X=\{\x\subset\R_+,\ |\x|<\infty\}$, where $|\x|$ means a number of points in $\x$.

As was said in the section \ref{s2}  the random events in the dynamics are separated  by a period of the deterministic motion during which configuration $\x$ moves into the positive direction as a rigid rod. The deterministic motion is described by the relations \reff{force} -- \reff{vel}. Recall that in the elastic case
\begin{equation}\label{vpath}
v(t)=v(0)e^{-\gamma\varkappa nt},
\end{equation}
where $n=|\x|$ (see \reff{vel}). The formula \reff{vpath} defines the plate dynamics when there are  no any random events on the interval $[0,t]$. It means that the number the contact points (which is equal to the element number in the set $\x$) is not changed on this interval and equal to $n$.

Further, a point $x\in\x$ of any configuration $\x$ we will call also \textit{particle}.

\subsubsection*{The stochastic Markov dynamic} There exist two kinds of the perturbations of the smooth deterministic dynamics: by a birth of a new particle or by a death of an existent particle.

The born particles are always localized at $0\in\R_+$ at the moment of their appearance. Thus it does not create immediately a resistant force. The birth intensity we denoted by $c_b = c_b(v(t))$.
The intensity $c_u = c_u(x(t))$ of the death of any particle may depend on the size $x(t)$ of the particle deformation.

Remark that the velocity $v(t)$ in the stationary state is always positive except the case $c_u(x)\equiv 0$.
%(\pch{Generally speaking, it is a theorem}). 
Any displacement $x\in \x$ cannot exceed the value $\left(\frac F \varkappa\right)^{\frac1\alpha}$.

The dynamics of the process is described by the following infinitesimal operator formalism.
First, we describe a set of  Markov process states. 

\vspace{.3cm}
\noindent
\textit{The  configuration set (the set of states)} we define as
\begin{eqnarray}\label{conf1}
\mathcal{X}&=& \bigcup_{n=0}^{\infty} \left\{ \{n\}\times\left[0,\frac{F}{\kappa}\right]^n\right\} \\
&=&\left\{(n,x_1,...,x_n):\:n\in\mathbb{N},\  \x=(x_1,...,x_n)\in\left[0,\left(\frac{F}{\varkappa}\right)^{\frac1\alpha}\right]^n \right\}
\nonumber
\end{eqnarray}
Every state $(n,x_1,...,x_n)$ means that there are $n$ the contact points (the particles) and $\x=(x_1,...,x_n)$ describes the deformations of the contacts.

\vspace{.3cm}

\noindent
\textit{Infinitesimal generator}. Let $\mathbf{H}=\{f=(f_n)\}$ be a set of continuous functions on $\mathcal{X}$, i.e. every $f_n:\:\{n\}\times\left[0,\frac{F}{\varkappa}\right]^n\to \R$ is continuous. We shall omit the index $n$ if it does not lead to misunderstanding and write $f(n,x_1,...,x_n)$ instead $f_n(n,x_1,...,x_n)$.
The infinitesimal operator $L$ of the Markov process  defined on $\mathbf{H}$ is
\begin{eqnarray}
&Lf(n,\mathbf{x})=&v\sum_{i=1}^n\frac{\partial f}{\partial x_i} \label{infsim} \\ 
&& {} +c_b(v)[f(n+1,x_1,...,x_n,0)-f(n,x_1,...,x_n)] \nonumber \\
&& {} + \sum_{j=1}^nc_u({x}_j)[f(n-1,x_1,...,\wh{x}_j,...,x_n)-f(n,x_1,...,x_n)],\nonumber
\end{eqnarray}
where $\wh{x}_j$ means that the variable $x_j$ is not presented  in the list of variables, and we recall that $v(n,x_1,...,x_n)= \gamma[F-\varkappa \sum_{i=1}^n x_i]_+$.
The first term on the right of \reff{infsim} corresponds to the deterministic plate motion between the random events. The second term reflects the birth event and the third term reflects the death event.

The Markov process defined by the operator $L$ is a piece-wise deterministic process (see \cite{Dav}).

The scaling limit of this Markov process is a deterministic process from the section \ref{s3} describing by the system \reff{y-3}, \reff{y-4} that follows from the general theory (see \cite{EK}

\subsection{Solution  \reff{vel} of equation \reff{dyn2}}\label{s41}Finding the solution of \reff{dyn2} on the interval $[t_1,t_2]$ introduce $X(t)=\sum_{\omega\in\Omega}x_\omega$. It follows from the elastic version of \reff{dyn2} that 
\begin{equation}\label{eq23}
nv(t)=\frac{\ed X(t)}{\ed t}=n\gamma\left(F-\varkappa X(t)\right),
\end{equation}
where $n=|\Omega|$ is the number of the contacts. The general solution \reff{eq23} is 
\begin{equation}\label{presol}
X(t)=\frac F\varkappa+Ce^{-n\gamma\varkappa (t-t_1)},
\end{equation}
where $C$ must be  defined from $X(t_1)$, that is $C=X(t_1)-\frac F\varkappa$. However if only the velocity  $v(t_1)$ is known then $C=-\frac{v(t_1)}{\varkappa\gamma}$. Now \reff{vel} follows from \reff{presol}.

\subsection{A dynamical system}\label{4.3} 
The stationary point \reff{stat} can be of two types: stable focus or stable node.
Essential role in the following plays the combination of the parameters
\begin{equation}\label{order}
a=\gamma\kappa \bar F\frac{\bar c_b}{\bar c_u^2},
\end{equation}
which we call an {\it order parameter}.

\begin{figure}
\begin{center}
\includegraphics[scale=0.4]{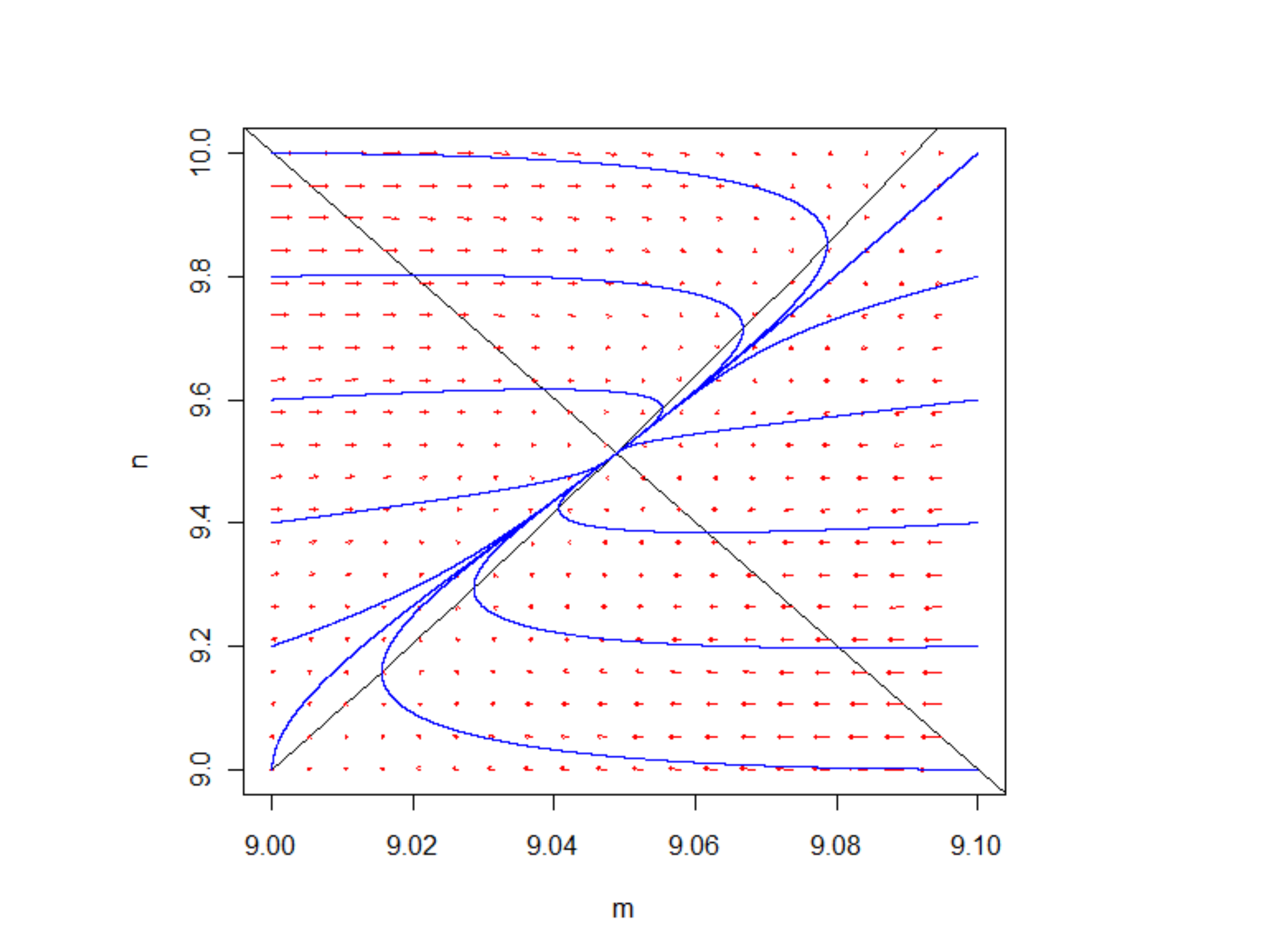}
\end{center}
\caption{The  integral curves of the field (\ref{y-7}), $a<20$.}%
\label{fig1}
\end{figure}
\begin{figure}
\begin{center}
\includegraphics[scale=0.4]{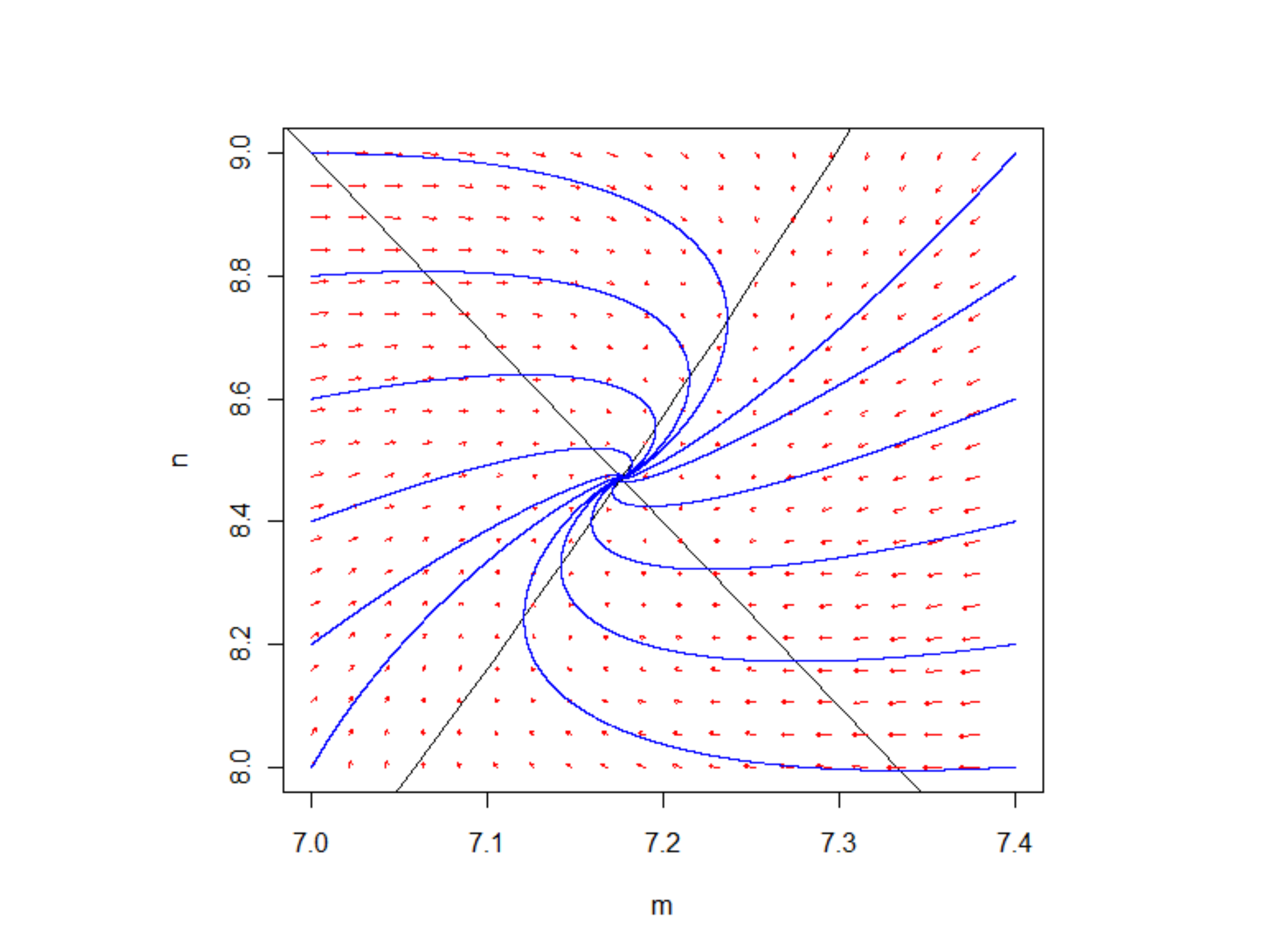}
\end{center}
\caption{The  integral curves of the field (\ref{y-7}), the case $a>20$.}%
\label{fig2}
\end{figure}

\begin{theorem}
There exists  in $\R^2_+$ an unique solution $(M_0,N_0)$ of the equations $dM/dt=dN/dt=0$ (\reff{y-7}):
\begin{eqnarray}\label{stat1}
% \nonumber to remove numbering (before each equation)
  M_0 &=& \frac{\bar F}{\kappa}- \frac{\bar F}{2\gamma\varkappa a}\left[\sqrt{1+4\gamma a}-1\right]  \\
  N_0 &=& \frac{\bar c_b}{\bar c_u}\frac{\bar F}{2a} \left[\sqrt{1+4\gamma a}-1\right]
 \nonumber
\end{eqnarray}
The linearized at $(M_0,N_0)$ equations \reff{y-7} are
\begin{eqnarray}\label{linear}
% \nonumber to remove numbering (before each equation)
  dM/dt &=&-(\bar c_u+\varkappa N_0 )(M-M_0) + (\bar F - \varkappa M_0)(N-N_0) \\
  dN/dt &=&  -\varkappa \bar c_b(M-M_0)-\bar c_u(N-N_0)
 \nonumber
\end{eqnarray}
The determinant of the matrix
$$
\mathcal{M}=\left(\begin{matrix}-(\bar c_u+\varkappa N_0 )& \bar F - \varkappa M_0 \\
-\varkappa \bar c_b& -\bar c_u \end{matrix}\right)
$$
is positive. Therefore the trace (the sum of the eigenvalues) of $M$ is negative.

If the order parameter $a<20$ then the eigenvalues are complex and the point $(M_0,N_0)$ is the stable focus, if $a>20$ the both eigenvalues are negative then the point $(M_0,N_0)$ is the stable node. The constant $a$ is defined by \reff{order} (see Figure \ref{fig1} and \ref{fig2}).
\end{theorem}
\section{Acknowledgement} E.P. and S.P. was partially supported by  RFBR Foundation (grant 13-01-12410). E.P. was supported by FAPESP (grant  2013/04040-7). A.Y. thanks CNPq (grant 307110/2013-3).

\end{document}